\documentstyle[12pt,amscd,amssymb]{amsart}
\newcommand{\rank}{\operatorname{rank}}

\newcommand{\codim}{\operatorname{codim}}
\newcommand{\pardeg}{\operatorname{pardeg}}

\newcommand{\Ker}{\operatorname{Ker}}
\newcommand{\fhom}{\operatorname{{\frak H}{\frak o}{\frak m}}}

\newcommand{\parhom}{\operatorname{{\frak P}{\frak a}{\frak r}{\frak H}{\frak
o}{\frak m}}}

\newcommand{\E}{{\cal E}}
\newcommand{\I}{{\cal I}}

\newcommand{\M}{{\cal M}}
\newcommand{\U}{{\cal U}}
\newcommand{\SU}{\operatorname{SU}}
\newcommand{\Aut}{\operatorname{Aut}}
\newcommand{\GL}{\operatorname{GL}}

\newcommand{\id}{\hbox{id}}
\newcommand{\al}{\alpha}
\newcommand{\be}{\beta}
\newcommand{\ga}{\gamma}

\newcommand{\ep}{\epsilon}

\newtheorem{defn}{Definition}[section]
\newtheorem{lem}[defn]{Lemma}
\newtheorem{thm}[defn]{Theorem}
\newtheorem{prop}[defn]{Proposition}
\newtheorem{cor}[defn]{Corollary}
\newtheorem{conj}[defn]{Conjecture}
\newtheorem{claim}[defn]{Claim}

\begin{document}
\title{Rationality of Moduli Spaces of Parabolic Bundles}
\author{Hans U. Boden and K\^{o}ji Yokogawa}
\address{Department of Mathematics and Statistics\\
	McMaster University \\
	Hamilton, Ontario \\
	Canada L8S 4K1}
\email{boden@@icarus.mcmaster.ca}
\address{Department of Mathematics \\
	Faculty of Science\\
	Osaka University\\
	Toyonaka, Osaka 560, Japan}
\email{yokogawa@@math.sci.osaka-u.ac.jp}
\thanks{
Both authors are members of the VBAC group of Europroj.}
\begin{abstract}
The moduli space of parabolic bundles with fixed determinant
over a smooth curve of genus greater than one is proved to be rational whenever
one of the multiplicities
associated to the quasi-parabolic structure is equal to one.
It follows that if rank and degree are coprime,
the moduli space of vector bundles is stably rational, and the bound
obtained on the level is strong enough to
conclude rationality in many cases.
\end{abstract}
\maketitle

\section{Introduction}

Let $X$ be a smooth complex curve of genus $g \geq 2$, $L$ a line bundle of
degree $d$ over $X,$ and
$\M_{r,L}$ the moduli space of semistable bundles $E$ of rank $r$ with
determinant $L$.
\begin{conj} \label{conj}
$\M_{r,L}$ is {\it rational}, i.e. it is birational
to a projective space.
\end{conj}
\noindent
Despite many positive results \cite{newstead1}, this is still an open problem,
even for $(r,d)=1.$

In this paper, we study a closely related problem, namely the birational
classification of moduli
spaces of parabolic bundles over $X.$
These moduli spaces occur naturally as
\begin{enumerate}
\item[(i)] unitary representation spaces of Fuchsian groups
\cite{mehta-seshadri},
\item[(ii)] moduli spaces of Yang-Mills connections on X with an orbifold
metric \cite{boden1}, and
\item[(iii)] moduli spaces of certain semistable bundles over an elliptic
surface \cite{bauer}.
\end{enumerate}
The theory developed in \cite{bh} and extended here shows that their birational
type depends only on the {\em quasi-parabolic} structure (see Proposition
\ref{prop:qp}). The methods of \cite{newstead1} then prove, in many cases, that
these moduli spaces are rational. The weaker result, Theorem \ref{thm:rat1},
uses only
Newstead's theorem, while the stronger one, Theorem \ref{thm:rat2}, requires
an adaptation of his inductive argument.

Using the theory developed in \S 4, it then follows from
Theorem \ref{thm:rat2} that
$\M_{r,L} \times {\Bbb P}^{r-1}$ is rational
if $(r,d)=1$ (see Corollary \ref{cor}).
Stable rationality of the moduli spaces
had been proved in this case by Ballico
\cite{ballico}, and  our result is a strengthening of his.
For instance, a consequence is that one can
conclude Conjecture \ref{conj}
under the assumption that $(r,d)=1$ for most values of the
genus\footnote{Choosing $d' \equiv d \mod(r)$ with $0<d'<r,$  the hypothesis is
that either $(d',g)=1$ or $(r-d',g)=1.$}
(see Corollary \ref{cor:rat}).

A number of useful facts are established along the way. One key point is
Proposition \ref{prop:fine},
which gives a simple criterion for the existence of a universal bundle of
stable parabolic bundles.
We also extend the theory developed in \cite{bh} in several important ways
(Theorems \ref{thm:bh1}, \ref{thm:bh2}, and \ref{thm:bh3});
the first two are standard but necessary for our purposes and the third is
completely new.
Its proof requires the idea of shifting a parabolic sheaf (Definition
\ref{defn:shift}), which also provides
a framework for the Hecke correspondence (equation (\ref{eqn:hecke})).
All of these results play a crucial role in the proofs of Theorems
\ref{thm:rat1} and \ref{thm:rat2}.

A brief word about the organization of this paper:
\S 2 introduces the notation used in the following sections,
\S 3 discusses the existence of universal families,
\S 4 summarizes and extends the theory of \cite{bh},
\S 5 describes shifting and the Hecke correspondence,
and \S 6 contains the proofs of the main results and their corollaries.

Before we begin, we would like to acknowledge a certain debt to the work of
Newstead,
upon which a number of our arguments depend,  and
without which this paper would be inconceivable.

\section{Notation}
Let $X$ be a smooth curve of genus $g \geq 2$ and  $D$ a reduced divisor in
$X.$
If $E$ is a ${\Bbb C}^r$ bundle over $X,$ then a {\it parabolic structure} on
$E$ with respect to
$D$ is just a collection of weighted flags in the fibers of $E$ over each $p
\in D$ of the form
\begin{eqnarray} \label{defn:parbun1}
&E_p = F_1(p) \supset F_2(p) \supset \cdots \supset F_{\kappa_p}(p) \supset
0,&\\
&\, \;\; 0 \leq \al_1(p) < \al_2(p)< \cdots < \al_{\kappa_p}(p) < 1.&
\label{defn:parbun2}
\end{eqnarray}
Holomorphic bundles $E$ with parabolic structures are called {\it parabolic
bundles}, and we use the notation $E_*$ to indicate the bundle (or,
equivalently, locally-free sheaf) $E$ together with a choice of parabolic
structure.
A morphism $\phi:E_* \longrightarrow E'_*$ of parabolic bundles is a bundle map
satisfying
$\phi(F_i(p)) \subset F'_{j+1}(p)$ whenever $\al_i(p) > \al'_j(p)$ for all $p
\in D.$
We use the tensor product notation $H^0(E_*^\vee \otimes E'_*)$
for these morphisms, where $E^\vee_*$ denotes the dual parabolic bundle (cf.
\cite{yoko}).

A {\it quasi-parabolic} structure on $E$ is what is left after the weights are
forgotten,
it is determined topologically by its flag type $m,$ which specifies {\it
multiplicities}
$m(p) = (m_1(p), \ldots, m_{\kappa_p}(p))$ for each $p \in D$ defined by
$m_i(p) = \dim F_i(p) - \dim F_{i+1}(p).$

A subbundle $E'$ inherits a parabolic structure from one on $E$ in a canonical
way:
The flag in $E'_p$ is gotten by intersecting with the flag in $E_p$ and the
weights are determined by choosing
maximal weights such that the inclusion map from $E'$ to $E$ is parabolic
(p.\ 213, \cite{mehta-seshadri}). Parabolic structures on quotients have a
similar description (loc.\ cit.).

A parabolic bundle $E_*$ is called {\it stable}
if every proper holomorphic subbundle $E'$ satisfies
$\mu (E'_*) < \mu (E_*),$
where
$$\mu (E_*) = \pardeg E_* / r =
 \deg E / r + \sum_{p \in D} \sum_{i=1}^{\kappa_p} m_i(p) \al_i (p)/r.$$
The parabolic bundle $E_*$ is called {\it semistable} if $\mu (E'_*) \leq \mu
(E_*)$
for each subbundle $E'_*.$
The construction of the moduli space $\M_\al$ of semistable parabolic bundles,
as a normal, projective variety, is given in \cite{mehta-seshadri}. The
subspace $\M^s_\al$
of stable bundles is smooth and Zariski-open, in particular,
if every semistable bundle is stable, then $\M_\al$ is smooth.

Let $\Delta^r = \{ (a_1,\ldots,a_r) \mid 0 \leq a_1 \leq \cdots \leq a_r <1\}$
and define $W = \{ \al:D\longrightarrow \Delta^r \}.$
Points in $W$ determine both the weights and the multiplicities. Conversely,
given a weight $\al$ in the sense of (\ref{defn:parbun2}), the associated point
in $W$ is
gotten by repeating each $\al_i(p)$ according to its multiplicity $m_i(p)$.
We abuse notation slightly by referring to points in $W$ as weights.
This gives an obvious notion of when a weight is {\it compatible}
with a choice of multiplicities, and for a given $m,$ we define
the open face of weights compatible with $m$ to be
$$V_m = \{ \al \in W \mid \al_{i-1}(p)=\al_i(p)
\Leftrightarrow \sum_{k=1}^j m_k(p) < i \leq \sum_{k=1}^{j+1} m_k(p) \}.$$

A weight in the interior of $W$ specifies full flags at each $p \in D.$
For every other choice of $m$,
$V_m$ is contained in the boundary of $W.$
Now $W$ is a simplicial set, and the face relations give a natural ordering on
$\{V_m\}$
and we write $V_m > V_{m'}$ if $V_{m'}$ is a proper face contained in the
closure of $V_m.$
This agrees with the natural ordering on $m$ gotten by successive refinement.

Weights for which $\M_\al$ is not necessarily smooth
satisfy $\mu(E'_*) = \mu(E_*)$ for some proper subbundle $E'.$
Letting $E''$ be the quotient, then the short exact sequence of parabolic
bundles
$E'_* \stackrel{\iota}{\longrightarrow} E_* \stackrel{\pi}{\longrightarrow}
E''_*$
determines a partition of $(d,r,m)$ in the obvious way:
$(d',d''), (r',r'')$ and $(m',m'')$ are the degrees, ranks, and multiplicities
of $(E',E'').$
(We define $m'$ and $m''$ here slightly unconventionally, namely
\begin{eqnarray*}
m'_i(p) &=& \dim (F_i(p) \cap \iota (E'_p)) - \dim (F_{i+1}(p) \cap \iota
(E'_p)),\\
m''_i(p) &=& \dim (\pi (F_i(p)) \cap E''_p) - \dim (\pi (F_{i+1}(p)) \cap
E''_p),
\end{eqnarray*}
for $p \in D$ and $1 \leq i \leq \kappa_p$.)
Notice that $r',r'' > 0$ and $m'_i(p), m''_i(p) \geq 0.$
Write $\xi = (d',r',m').$ For fixed
$\xi,$ the set of weights
compatible with $m$
for which $\mu(E'_*) = \mu(E_*)$ is
the intersection of a hyperplane $H_\xi$ in $W$ with $V_m$ given by the
equation
\begin{equation}\label{eqn:hyper}
\sum_{i=1}^{\kappa_p} m_i(p) \al_i(p))
\sum_{p\in D} \sum_{i=1}^{\kappa_p} (r' m_i(p)-r m'_i(p)) \al_i(p) = r d' -
r'd.
\end{equation}
There are only finitely many hyperplanes; the above equation puts a bound
on $d'$ and all other quantities are already bounded.
We shall refer to $H_\xi \cap V_m$ as a {\it wall} in $V_m.$
These walls induce a chamber structure on $V_m,$ a {\it chamber} being
a connected component of $V_m \setminus \cup_\xi H_\xi$
(it is possible that $V_m \subset H_\xi$).
Weights $\al \in W \setminus \cup H_\xi$ are called {\it generic}, and for
these weights,
$\M_\al = \M^s_\al.$
In the next section, we shall see that $V_m$ contains a generic weight
if and only if the degree $d$ and the set of multiplicities $\{m_i(p)\}$
have  greatest common divisor equal to one.

\section{Families of parabolic bundles}
In this section, we present Proposition \ref{prop:fine}, which establishes
the existence of a universal family of stable parabolic bundles
parametrized by $\M_\al^s$ whenever $V_m$ contains a generic weight.
Although results of this type are well-known to experts,
the proposition, as well as the proof, are original
(cf. Th\'eor\`eme 32, \cite{seshadri}).
It is important because, in the case of ordinary bundles,
the non-existence of the universal family
(\cite{ramanan}) is the obstruction to
proving Corollary \ref{conj} by induction,
and as shown in \S 6,
the analogous argument works for parabolic
bundles precisely because the necessary conditions
for the vanishing of this obstruction given by Proposition \ref{prop:fine}
are often satisfied.

Given positive integers $m_1,\ldots,m_\kappa$ such that $m_1+\cdots m_\kappa =
r,$
define ${\cal F}_m$ to be the variety of flags of type $m.$
These are simply flags
${\Bbb C}^r=F_1 \supset \cdots \supset F_s \supset 0$ with $\dim F_i -\dim
F_{i+1} = m_i.$
Furthermore, for any bundle
$E \longrightarrow S$ of rank $r,$ let ${\cal F}_m(E) \longrightarrow S$ be the
bundle of flags of type $m.$

Given a bundle $U \rightarrow S \times X,$
we adopt the notation $U_s = U|_{\{s\} \times X}.$ We also use
$\pi_S$ for the projection map $S\times X \rightarrow S.$

\begin{defn} \label{defn:family}
Fix multiplicities $m(p)$ for each $p \in D.$
\begin{enumerate}
\item[(i)]
A family of quasi-parabolic bundles (of type $m$) parametrized by a variety
$S$ is a bundle $U$ over $S \times X$ together with a section $\phi_p$ of
the flag bundle ${\cal F}_{m(p)}(U|_{S \times \{p\}}) \longrightarrow S$ for
each $p \in D.$
\item[(ii)] Two families $(U,\phi)$ and $(U',\phi')$ parametrized by $S$ are
equivalent, written
$(U,\phi) \sim (U',\phi'),$ if there exists a line bundle $L$ over $S$ and an
isomorphism
$U \cong U' \otimes \pi_S^* L$ under which $\phi \mapsto \phi'.$
\end{enumerate}
\end{defn}

Note that the section $\phi_p$ in (i) above is just a choice of
a nested chain of subbundles of $U|_{S \times \{p\}}$ whose relative coranks
are given by the multiplicities $m(p).$
A family of parabolic bundles is gotten by associating a fixed set of
weights to each chain of subbundles.
Let $U_*=(U,\phi,\al)$ be the resulting family of parabolic bundles
and $U_{s,*}=(U_s, \phi(s),\al)$ be the parabolic bundle above $s \in S.$
Then $U_*$ is called a family of (semi)stable parabolic bundles if
$U_{s,*}$ is (semi)stable for each $s \in S.$

It follows from the construction of Mehta and Seshadri that $\M_\al$ is a
coarse moduli space.
Proposition 1.8 of  \cite{newstead2} then gives two conditions which are
necessary and sufficient  for a coarse moduli space to be fine, i.e. to admit a
universal family.
The second condition is not difficult to verify
using an argument similar to that given in Lemma 5.10 of \cite{newstead2}.
The first condition requires that we construct a family
$\U^\al_*$ parametrized by $\M^s_\al$ with the property that $\U^\al_{e,*}$
is a parabolic stable bundle isomorphic to $E_*$ for all
$[E_*] = e \in \M^s_\al.$

To construct this family, we need to review the construction of
$\M_\al$ (\cite{bhosle}, \cite{mehta-seshadri}).
Let $Q$ be the Hilbert scheme of coherent sheaves over $X$ which are quotients
of ${\cal O}_X^{\oplus N}$ with fixed Hilbert polynomial
(that of $E(k)$ for $k \gg g$), where $N= h^0(E).$ Let $U$ be the universal
family on
$Q \times X.$ Define $R$ to be the subscheme of $Q$ of points $r \in Q$ so that
$U_r$ is a locally free sheaf which is generated by its global sections and
$h^1(U_r)=0.$
Let $\widetilde{R}$ be the total space of the universal flag
bundle over $R$ with flag type $\prod_{p\in D} {\cal F}_{m(p)},$
and let $\widetilde{U}$ be the pullback of $U$ to $\widetilde{R}.$
Then $\widetilde{U}$ is canonically a family of parabolic bundles
parametrized by $\widetilde{R}$ by letting, for each $p \in D,$
$\phi_p$ be the tautological section
and $\al(p)$ be the fixed weights.
It follows that $\widetilde{R}$ has the local universal property for parabolic
bundles
(p.\ 16, \cite{bhosle}).

The subsets $\widetilde{R}^{s}$ ($\widetilde{R}^{ss}$) corresponding to the
stable (semistable)
parabolic bundles are invariant under the natural action of $\GL(N) =
\Aut({\cal O}_X^{\oplus N}),$
and $\M_\al$ is a good quotient of $\widetilde{R}^{ss}$
(with linearization induced by the weights $\al$),
and $\M^s_\al$ is the geometric quotient of
$\widetilde{R}^{s}.$

The center of $\GL(N)$
acts trivially on $R$ and $\widetilde{R},$ but nontrivially on the locally
universal bundle $\widetilde{U}$.
In fact, $\lambda(\id)$ acts on $\widetilde{U}$ by scalar multiplication by
$\lambda$
in the fibers (this follows from p.\ 138, \cite{newstead2}).
Given a line bundle $L$ over $\widetilde{R}^{s}$ with a natural
lift of the $\GL(N)$ action such that $\lambda(\id)$ acts by multiplication by
$\lambda,$
then using $\widetilde{U}^{s}$ to denote $\widetilde{U}|_{\widetilde{R}^{s}
\times X},$
the quotient of $\widetilde{U}^{s} \otimes \pi^*_{\widetilde{R}^{s}} L^{-1},$
together with the tautological sections and weights $\{\phi_p, \al(p) \mid p
\in D \}$ mentioned above,
gives the desired family.

\begin{prop} \label{prop:fine}
Such a line bundle $L$ exists if either
\begin{enumerate}
\item[(i)] the elements of the set $\{d, m_i(p) \mid p \in D, 1 \leq i \leq
\kappa_p \}$
have greatest common divisor equal to one, or
\item[(ii)] the face $V_m$ containing $\al$ contains a generic weight.
\end{enumerate}
Moreover, these two conditions are equivalent, and when they are satisfied,
the moduli space $\M^s_\al$ is fine.
\end{prop}

The idea of the proof is to find line bundles $L_k$ for each $k \in \{d, m_i(p)
\}$
over $\widetilde{R}^{s}$ with natural actions of $\GL(N)$ such that
$\lambda(\id)$ acts
by scalar multiplication by $\lambda^k.$
Then (i) gives the existence of $k_1, \ldots, k_\ell \in \{d, m_i(p) \}$
and integers $a_1, \ldots, a_\ell$ so that
$a_1 k_1 + \cdots a_\ell k_\ell = 1.$
The required line bundle is then the tensor product
$L =  L^{a_1}_{k_1} \otimes \cdots \otimes L^{a_\ell}_{k_\ell}.$
At the end of the proof, we will show that (i) and (ii) are equivalent.

We start with a lemma.
\begin{lem}
Suppose $E_*$ is parabolic semistable or degree $d$ and rank $r$ and $H_*$ is a
parabolic line
bundle of degree $h$, then
\begin{equation} \label{lem:ineq}
h^1(H^\vee_* \otimes E_*) \neq 0 \quad \Rightarrow \quad d \leq r(2g-2+h) + r^2
n.
\end{equation}
\end{lem}
\begin{pf}
Serre duality for parabolic bundles (Proposition 3.7 of \cite{yoko}) implies
that
$$h^1(H^\vee_* \otimes E_*)
 \leq h^0(E^\vee_* \otimes {H}_* \otimes K(D)).$$
(If we had used $h^0(E^\vee_* \otimes \widehat{H}_* \otimes K(D)),$
the circumflex over $H_*$ indicating {\it strongly} parabolic morphisms, we
would get the usual statement
of Serre duality with
equality, cf. \cite{yoko,by}.)
Suppose that $\phi : E \longrightarrow H \otimes K(D)$ is a non-zero map and
let $E'$ be the subbundle generated by $\Ker \phi.$ Then
$$\deg E' \geq \deg E - \deg H\otimes K(D) = d - h - (2g -2 +n).$$
Considering $E'_*$ with its canonical parabolic structure as a subbundle of
rank $r-1,$
the inequality (\ref{lem:ineq}) follows easily from this,
semistability of $E_*,$ and
the inequalities $\pardeg E'_* \geq \deg E'$ and
$\pardeg E_* \geq \deg E + rn.$
\end{pf}

\noindent
{\it Proof of Proposition.}
Write the weights $\al$ without repetition.
Choose $\ell: D \longrightarrow {\Bbb Z}$
with
$ 1 \leq \ell_p \leq \kappa_p+1$ and set $\be(p) = \al_{\ell_p}(p).$
(Take $\be(p) > \al_{\kappa_p}$ if $\ell_p = \kappa_p+1.$)
For $h \in {\Bbb Z},$ define
$$\chi(\ell,h) = d + r(1-g-h) - \sum_{p \in D} \sum_{i=1}^{\ell_p-1} m_i(p).$$
Let $H_*$ be the parabolic line bundle with $\deg H = h < d/r - r n -(2g-2)$
and with weights
$\be(p)$ at $p \in D.$ It follows from the lemma that
if $E_*$ is semistable, then $h^1(H^\vee_* \otimes E_*)=0.$ Thus
$h^0(H^\vee_* \otimes E_*) = \chi(\ell,h)$ by Riemann-Roch.
Hence $(R^0 \pi_{\widetilde{R}^{s}}) (\widetilde{U}^{s} \otimes \pi_X^* H_*)$
is a locally free sheaf of rank $\chi(\ell,h)$ over $\widetilde{R}^{s}.$
Let $L(\ell,h)$ be the determinant of the corresponding bundle.
By construction, the $\GL(N)$ action on $\widetilde{U}$ induces one on this
bundle (and hence on $L(\ell,h)$);
$\lambda(\id)$ acts by scalar multiplication by $\lambda$ on the bundle and by
$\lambda^{\chi(\ell,h)}$ on $L(\ell,h).$
It is now a simple exercise in high school algebra to see that we can choose
$h,h'$ and $\ell,\ell'$ so
that $\lambda(\id)$ acts on $L(\ell,h) \otimes L(\ell',h')$ by $\lambda^k$ for
any $k \in \{ d, m_i(p) \}.$

This proves the conclusion of the proposition assuming (i), and now we show
that conditions (i) and (ii)
are equivalent. Suppose first that (i) does not hold.
Consider $E_*$ as
a quasi-parabolic bundle without holomorphic structure, which will be specified
later.
Since the set $\{ d, m_i(p)\}$ is not relatively prime,
there exists a prime number $q$
evenly dividing each element of the set. Clearly $q$ also divides
$r.$
Set $d' =d/q, r' = r/q$ and $m'_i(p) =m_i(p)/q.$
Consider the quasi-parabolic bundle $E'_*$
with degree $d'$, rank $r',$  and multiplicities $m'.$
Any choice of weights $\al$ on $E_*$ induces (the same!) weights on $E'_*,$ and
it follows that since $g \geq 2,$ there is some holomorphic structure for
which $E'_*$
is semistable.
Define the holomorphic structure on $E_*$ by
$$E_* = E'_* \oplus \stackrel{q}{\cdots} \oplus E'_*.$$
It follows that  $E_*$ is semistable but not stable
for {\it any} choice of compatible weights.
This implies that $V_m$ does not contain a generic weight.

Suppose conversely that $V_m$ does not contain a generic weight.
Since $V_m$ is affine,
$V_m \subset H_\xi$ for some $\xi=(r',d',m')$ Using (\ref{eqn:hyper}), we
conclude that
for all $\al \in V_m,$
$$ \sum_{p \in D} \sum_{i=1}^{\kappa_p} (r m_i'(p) - r' m_i(p)) \al_i(p) =
rd'-r'd .$$
(Here, we are still thinking of $\al$ without repetition.)
We can vary each $\al_i(p)$ continuously by some small amount,
and it follows that $$rm'_i(p) - r' m_i(p)=0= rd'-r'd$$
for all $i$ and $p.$
Since $r' < r,$ there exists a prime $q$ such that
$q^k$ divides $r$ but not $r'.$ Hence $q$ divides $d$ and each element of
the set $\{m_i(p) \mid p \in D, 1 \leq i \leq \kappa_p\}.$
$\quad \Box$
\section{The variation and degeneration theorems}

In this section, we describe and extend the theory of \cite{bh}. This allows us
to
compare the moduli spaces of parabolic bundles $\M_\al$ and $\M_\be$ when
\begin{enumerate}
\item[(i)] $\al, \be \in V_m$ are generic weights in adjacent chambers,
\item[(ii)] $\al \in V_\ell$ and $\be \in V_m$ are generic weights not
separated by any hyperplanes and
$V_\ell > V_m.$
\end{enumerate}
Cases (i) and (ii) correspond to Theorem 3.1 and Proposition 3.4 of \cite{bh}.
We present slightly stronger versions of those results tailored for our
purposes here.

Starting with (i), suppose that $\al, \be \in V_m$ are generic weights
separated by a single hyperplane $H_\xi.$
Choose $\ga \in H_\xi$ on the straight line connecting $\al$ to $\be.$
Then $\M_\ga$
is stratified by the Jordan-H\"older type of the underlying bundle,
and since $\ga$ lies on only one hyperplane, there are exactly two strata:
the stable bundles $\M^s_\ga$ and the strictly semistable bundles
$\Sigma_\ga.$
Writing $\xi=(r',d',m')$
for the partition, then it is not hard to see that
$\Sigma_\ga \cong \M_{\ga'} \times \M_{\ga''},$ with the obvious definitions
for $\ga'$ and $\ga''$
coming from the partition $\xi.$
\begin{thm} \label{thm:bh1}
There are natural algebraic maps $\phi_\al$ and $\phi_\be$
$$\begin{array}{rcl}\M_\al& & \M_\be\\
\;\; \phi_\al \!\! \searrow \!\!\!\!\!\!\! && \!\!\!\!\!\!\! \swarrow \!\!
\phi_\be\\& \M_\gamma
\end{array}$$
which are generized blow-downs along projectivizations of vector bundles over
$\Sigma_\ga,$
where the projective fiber dimensions $e_\al$ and $e_\be$ satisfy $e_\al +
e_\be + 1 = \codim \Sigma_\ga.$
\end{thm}
\begin{pf}
The proof is the same as in \cite{bh}, the only difference being the actual
computation of the
numbers $e_\al$ and $e_\be,$ which we discuss now.
We assume that $E_* \sim_S E'_* \oplus E''_*,$ where $[E_*] \in \Sigma_\ga$ and
$\sim_S$ denotes Seshadri equivalence (i.e.\ isomorphic Jordan-H\"older form).
The topological type of the parabolic bundles $E'_*$ and $E''_*$ does not
change as
$[E_*]$ varies within $\Sigma_\ga.$ We use $(r',r''), (d',d'')$ and $(m',m'')$
to denote
the ranks, degrees, and multiplicities of $(E'_*, E''_*),$
written as in \S 2.
The moduli spaces $\M_\al, \M_\be,$ and $\M_\ga$ have dimension
$$(g-1) r^2 +1 + \frac{1}{2}\sum_{p \in D} \left(r^2 - \sum_{i=1}^{\kappa_p}
m_i(p)^2\right).$$
Using a similar formula for $\Sigma_\ga=\M^{\ga'} \times \M^{\ga''},$ we find
that
$$\codim \Sigma_\ga = r' r''(2g - 1) -1+\sum_{p\in D}\sum_{i=1}^{\kappa_p}
m'_i(p) m''_i(p).$$
Now we claim that
$$h^0({E''_*}^\vee \otimes E'_*)=0=h^0({E'_*}^\vee \otimes E''_*).$$
This is true for any $\al' \in V_m,$ as one of these equations is true for
$\al,$ the other for $\be,$
but $H^0$ is constant as the weights are varied within $V_m.$
Let $\U'$ and $\U''$ be the families parametrized by $\Sigma_\ga$ gotten by
pulling back the universal
families $\U^{\ga'}$ and $\U^{\ga''},$ whose existence follows from Proposition
\ref{prop:fine}.
Then the vector bundles referred to in the theorem are
$$(R^1 \pi_{\Sigma_\ga})({\U''}^\vee \otimes \U') \hbox{ and } (R^1
\pi_{\Sigma_\ga})({\U'}^\vee \otimes \U'').$$
The projectivizations of these  bundles have dimensions
\begin{eqnarray}
e_\al &=& h^1({E''_*}^\vee \otimes E'_*)-1 = r'' d' - r' d'' + r' r'' (g-1) +
\chi ({\cal Q}) -1,
\label{formula:ealpha}\\
e_\be &=& h^1({E'_*}^\vee \otimes E''_*)-1 = r' d'' - r'' d' + r' r'' (g-1) +
\chi ({\cal Q}') -1,
\label{formula:ebeta}
\end{eqnarray}
where
${\cal Q}$ and ${\cal Q}'$ are skyscraper sheaves supported on $D$
obtained as the quotients
\begin{eqnarray*}
&\parhom (E''_*,E'_*) \longrightarrow \fhom (E'',E') \longrightarrow {\cal
Q},&\\
&\parhom (E'_*,E''_*) \longrightarrow \fhom (E',E'') \longrightarrow {\cal
Q}'.&
\end{eqnarray*}
It is a nice exercise to see
$$ \chi ({\cal Q}) + \chi ({\cal Q}') =  \sum_{p \in D}
\left( r' r'' -\sum_{(i,j) \in S_e(p)} m'_i(p)m''_j(p) \right),$$
where $S_e(p) = \{ (i,j) \mid \ga'_i(p) = \ga''_j(p)\}.$
This shows $e_\al + e_\be +1 = \codim \Sigma_\ga.$
\end{pf}

\begin{thm} \label{thm:bh2}
Suppose that $\al \in V_\ell, \, \be \in V_m, \, V_\ell > V_m,$ and that
$\al$ and $\be$ are generic and are not separated by any hyperplanes.
Then there exists a fibration $\psi:\M_\al \longrightarrow \M_\be$
with fiber a (possibly twisted) product of flag varieties and
this fibration is locally trivial in the Zariski topology.
In particular, $\M_\al$ is birational to the product of $\M_\be$ with a
product of flag varieties.
\end{thm}
\begin{pf}
The hypothesis $V_\ell > V_m$ just means that the flag
structure degenerates as we pass from $\al$ to $\be.$
By induction, it is enough to prove the above statement
when the degeneration of the flag structure
is taking place at only one parabolic point.
Given
$E_*$ a parabolic bundle with multiplicities $m$ and weights $\al,$
let $E'_*$ be the parabolic bundle
with multiplicities $\ell$ and weights $\be$
resulting from forgetting part of the flag structure and interchanging the
weights.
One easily verifies that if $E_*$ is $\al$-stable, then $E_*'$ is $\be$-stable,
and the existence of the morphism $\psi$ then follows from the coarseness
property of $\M_\be.$

The remaining issue is to identify the fiber
and to prove local triviality.
For the first issue, notice that there is an inverse procedure to the forgetful
map described above.
Given a parabolic bundle $E'_*$ with multiplicities $\ell$  and weights $\be,$
consider all parabolic bundles $E_*$ with weights $\al$ obtained  from $E'_*$
by
refining the flag stucture to one with multiplicities $m$ and
exchanging the weights.
For a given $E'_*,$ the set of all such possible refinements $E_*$ is
parametrized by a flag variety.

A straightforward numerical verification shows that applying this procedure
to a $\be$-stable parabolic bundle $E_*'$ yields
an $\al$-stable $E_*$ for every possible refinement.
It is not hard to see that the same procedure, when applied to the universal
family $\U_*^\be,$
identifies $\M_\al$ with the total space of the flag bundle of
$\U^\be$ restricted to $\M_\be \times \{p\}$ and the map $\psi$
with the bundle projection.
\end{pf}


One might expect from Theorem \ref{thm:bh1} that the birational type of
$\M_\al$ depends only on the underlying quasi-parabolic structure.
This is the content of the following proposition.

\begin{prop} \label{prop:qp}
Suppose that $g \geq 2.$
Then the birational type of $\M_\al$ is independent of the
choice of $\al \in V_m.$
\end{prop}

\begin{pf}
We prove the proposition by showing that $\M_\al$ and $\M_\be$
are birational whenever
$\al,\be \in V_m$
are not separated by any walls
(although one may lie on a wall which does not contain the other).
So assume that
$\al \in \cap_{i=1}^n H_{\xi_i}$ and $\be \in \cap_{i=1}^m H_{\xi_i},$
where $m \geq n.$
By Theorem 4.1 \cite{mehta-seshadri}, $\M_\al$ and $\M_\be$ are normal,
projective varieties and
$\dim \M_\al = \dim \M_\be,$ hence we only need to construct an
injective morphism $\phi:\M^s_\be \longrightarrow \M^s_\al$
to conclude $\M_\al$ is birational to $\M_\be.$
One easily verifies that every $\be$-stable bundle is $\al$-stable,
and the existence of $\phi$ follows from the
coarseness of $\M_\al.$
\end{pf}

\section{Shifting and the Hecke correspondence}
In this section, we introduce the notion of a shifted parabolic bundle,
which is the result of changing the weights, multiplicities, and degree of
$E_*$
in a prescribed way.
In some sense, shifting is a symmetry of a larger weight space, one
which includes bundles of different degrees. Two applications of shifting are
discussed at the end.

Shifting is most naturally described in terms of parabolic sheaves.
If $\E$ is a locally free sheaf on $X,$
then a {\it parabolic structure} on $\E$ consists of a weighted filtration of
the form
\begin{eqnarray} \label{eqn:filtration1}
&\E=\E_{\al_1} \supset\E_{\al_2} \supset \cdots \supset \E_{\al_l} \supset
\E_{\al_{l+1}} = \E(-D),&\\
&0 \leq \al_1 < \al_2 < \cdots < \al_l < \al_{l+1}=1.\quad \quad&
\label{eqn:filtration2}
\end{eqnarray}
We can define $\E_x$ for $x \in [0,1]$ by setting $\E_x = \E_{\al_i}$
if $\al_{i-1} < x \leq \al_i,$ and then extend to $x \in {\Bbb R}$ by setting
$\E_{x+1} =\E_{x}(-D).$
We call the resulting filtered sheaf $\E_*$ a parabolic sheaf and $\E=\E_0$
the underlying sheaf.

We can define parabolic subsheaves, degree, and stability for
these objects, and there
is a categorical equivalence between locally free parabolic sheaves and
parabolic bundles.
We describe this in case $D = p,$ the general case being quite similar
(\cite{yoko}, \cite{by}).

Suppose that $E_*$ is a parabolic bundle given by flags and weights in the
fibers as in
(\ref{defn:parbun1}) and (\ref{defn:parbun2}).
Define $\E_*$ by setting $$\E_x = \ker(E\rightarrow E_{p}/F_i),$$ for $
\al_{i-1} < x < \al_i.$
Thus $\E_*$ is a parabolic sheaf.
Conversely, given a parabolic sheaf $\E_*,$ the quotient $\E_0/\E_1 =
\E/\E(-p)$
is a skyscraper sheaf with support $p$ and fiber that of $\E.$
Defining  a flag in this fiber by setting
$F_i = (\E_{\al_{i}}/\E_1)_p$ and associating the weight $\al_i,$
we obtain a parabolic bundle in the sense of (\ref{defn:parbun1}) and
(\ref{defn:parbun2}).

The category of parabolic sheaves is developed in \cite{yoko}, where one finds
for example
the definitions of tensor products $\E_* \otimes \E'_*$ and duals $\E^\vee_*.$
We use this notation freely in the
calculations of \S 6 involving sheaf cohomology and point out that
$H^i(\E_*) = H^i(\E).$

\begin{defn} \label{defn:shift}
Given a parabolic sheaf $\E_*$ and $\eta \in {\Bbb R},$ define the shifted
parabolic sheaf
$\E_*[\eta]_*$ by setting $\E_*[\eta]_x = \E_{x+\eta}.$
\end{defn}
\noindent
{\it Remark.} The above operation can be refined in case $D = p_1 + \cdots +
p_n.$ If
$\eta=(\eta_1, \ldots, \eta_n),$ then one can shift $\E_*$ by $\eta_i$ at each
$p_i \in D$
(\cite{yoko}, \cite{bh}).

\medskip \noindent
It is not difficult to verify that $\E_*[\eta]_*$ is (semi)stable if and only
if $\E_*$ is (semi)stable,
and it follows that this defines an isomorphism between the associated moduli
spaces of parabolic bundles.

We can easily describe the parabolic structure on
the shifted bundle $\E'_* = \E_*[\eta]_*$ in case $0 < \eta \leq 1$ and $D =
p.$
Let $E'_*$ denote the parabolic bundle associated to $\E'_*.$
If $i$ is the integer with $\al_i < \eta \leq \al_{i+1},$ then
the weights of $E'_*$ are given by
\begin{equation} \label{eqn:sft_wts}
\al'_j = \begin{cases} \al_{j+i} - \eta & \hbox{ for  } j=1, \ldots, r-i, \\
                        1+\al_{j-r+i}-\eta & \hbox{ for  }  j = r-i+1,\ldots,
r.
\end{cases}
\end{equation}
The quasi-parabolic structure of $E'_*$ has
multiplicities $m'$ given by a cyclic permutation of $m,$
i.e.\ $m' = (m_{i+1}, \ldots, m_\kappa, m_1, \ldots, m_i).$
Although $\E'$ is a subsheaf of $\E,$
$E'$ is {\it not} a subbundle of $E,$ so one must appeal to sheaf
theory in order to define the flag in
$E'_p.$ This is a simple exercise in tracing
through the equivalence between locally free parabolic sheafs and
parabolic bundles given above.

\begin{figure}[b]
\begin{picture}(100,210)(-110,-10)
        \put(-70,150){$\E_*$}
        \put(-70,100){${\Bbb R}$}
        \put(-40,110){\line(1,0){290}}     
        \put(5,110){\line(0,1){70}}      		
        \put(3,100){$0$}
        \put(170,110){\line(0,1){20}}      		
        \put(168,100){$1$}
        \put(-12,185){$\E=\E_{\al_1}$}
        \put(-40,180){\line(1,0){78}}      
        \put(40,180){\circle*{4}}        
        \put(34,100){$\al_1$}
        \multiput(40,110)(0,4){12}{\line(0,1){2}}        
        \multiput(40,178)(0,-4){4}{\line(0,-1){2}}       
        \put(55,165){$\E_{\al_2}$}
        \put(42,160){\line(1,0){48}}     
        \put(40,160){\circle{4}}         
        \put(90,160){\circle*{4}}        
        \put(88,100){$\al_2$}
        \multiput(90,110)(0,4){7}{\line(0,1){2}}         
        \multiput(90,158)(0,-4){4}{\line(0,-1){2}}       
        \put(92,140){\line(1,0){48}}     
        \put(109,145){$\E_{\al_3}$}
        \put(90,140){\circle{4}}         
        \put(140,140){\circle*{4}}       
        \put(137,100){$\al_3$}
        \multiput(140,110)(0,4){5}{\line(0,1){2}}        
        \multiput(140,138)(0,-4){2}{\line(0,-1){2}}      
        \put(142,130){\line(1,0){48}}    
        \put(140,130){\circle{4}}        
        \put(190,130){\circle*{4}}       
        \put(152,135){$\E(-p)$}
        \put(185,100){$1\!+\!\al_1$}
        \multiput(190,110)(0,4){2}{\line(0,1){2}}
        \multiput(190,128)(0,-4){2}{\line(0,-1){2}}
        \put(192,120){\line(1,0){58}}
        \put(190,120){\circle{4}}
        \put(200,125){$\E_{\al_2}(-p)$}
        \put(-80,50){$\E_*[\eta]_*$}
        \put(-70,0){${\Bbb R}$}
        \put(-40,10){\line(1,0){250}}     
        \put(5,10){\line(0,1){50}}      		
        \put(3,0){$0$}
        \put(170,10){\line(0,1){10}}      		
        \put(168,0){$1$}
        \put(-40,85){$\E_{\al_1}$}
        \put(-40,80){\line(1,0){18}}      
        \put(-20,80){\circle*{4}}        
        \put(-40,0){$\al_1\!-\!\eta$}
        \multiput(-20,10)(0,4){12}{\line(0,1){2}}        
        \multiput(-20,78)(0,-4){4}{\line(0,-1){2}}       
        \put(-5,65){$\E_{\al_2}$}
        \put(-18,60){\line(1,0){48}}     
        \put(-20,60){\circle{4}}         
        \put(30,60){\circle*{4}}        
        \put(20,0){$\al_2\!-\!\eta$}
        \multiput(30,10)(0,4){7}{\line(0,1){2}}         
        \multiput(30,58)(0,-4){4}{\line(0,-1){2}}       
        \put(32,40){\line(1,0){48}}     
        \put(49,45){$\E_{\al_3}$}
        \put(30,40){\circle{4}}         
        \put(80,40){\circle*{4}}       
        \put(67,0){$\al_3\!-\!\eta$}
        \multiput(80,10)(0,4){5}{\line(0,1){2}}        
        \multiput(80,38)(0,-4){2}{\line(0,-1){2}}      
        \put(82,30){\line(1,0){48}}    
        \put(80,30){\circle{4}}        
        \put(130,30){\circle*{4}}       
        \put(92,35){$\E(-p)$}
        \put(110,0){$1\!+\!\al_1\!-\!\eta$}
        \multiput(130,10)(0,4){2}{\line(0,1){2}}
        \multiput(130,28)(0,-4){2}{\line(0,-1){2}}
        \put(132,20){\line(1,0){78}}
        \put(130,20){\circle{4}}
        \put(140,25){$\E_{\al_2}(-p)$}

\end{picture}
\caption{The parabolic sheaf $\E_*$ shifted by $\eta$ with $\al_1 < \eta <
\al_2.$}
\end{figure}
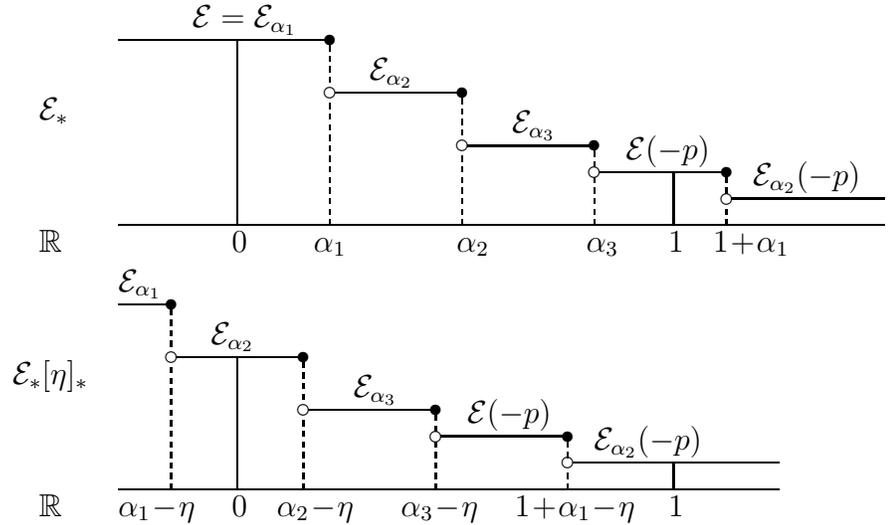

There are two interesting applications of shifting we discuss now.
The first is the Hecke correspondence. Using $\M_{r,d}$ to denote the moduli
space of
semistable bundles of rank $r$ and degree $d,$ the Hecke correspondence gives a
means
of comparing $\M_{r,d}$ and
$\M_{r,d'}$ through
the use of parabolic bundles. For $r=2,$ this was observed in a remark  at the
end of \cite{mehta-seshadri}.

To start, define $\ep_+(d,r), \ep_-(d,r),$ and $\ep(d,r)$  for $d,r \in {\Bbb
Z}$ with $r >0$ by
\begin{eqnarray*}
 \ep_\pm(d,r) &=& \inf \{ \pm({\textstyle \frac{d}{r}-\frac{d'}{r'}}) \mid
d',r' \in {\Bbb Z}, \; 1 \leq r' < r, \hbox{ and } \pm ({\textstyle
\frac{d}{r}-\frac{d'}{r'}}) > 0 \} \\
 \ep(d,r) \;\; &=& \min \{ \ep_\pm(d,k) \mid k=1, \ldots, r \}.
\end{eqnarray*}
It is easy to see that $\ep_\pm(d,k) > 0$ for all $k,$ thus $\ep(d,r)>0$ as
well.

Suppose that $E$ is a bundle over $X$ of degree $d$ and rank $r$
and suppose further that $E'$ is a proper subbundle.
If $\mu(E') < \mu(E),$ then $\mu(E)-\mu(E') \geq \ep_+(d,r).$
Similarly, if $\mu(E') > \mu(E),$ then $\mu(E')-\mu(E) \geq \ep_-(d,r).$

\begin{prop} \label{prop:num}
Suppose that $E_*$ satisfies $ {\displaystyle \sum_{p \in D}
\sum_{i=1}^{\kappa_p} m_i(p)} \al_i(p) < \ep(d,r)/2.$
\begin{enumerate}
\item[(i)] If $E$ is stable as a regular bundle, then $E_*$ is parabolic
stable.
\item[(ii)] If $E_*$ is parabolic stable, then $E$ is semistable as a regular
bundle.
\end{enumerate}
\end{prop}
\begin{pf}
(i)
If $E'_*$ is a proper parabolic subbundle of $E_*$, then
$$\mu(E'_*) \leq \mu(E') + \ep(d,r)/2 < \mu(E') + \ep_+(d,r) \leq \mu(E) <
\mu(E_*),$$
thus $E_*$ is parabolic stable.

(ii) If $E'$ is a subbundle of $E,$ then
$$\mu(E') \leq \mu(E'_*)  < \mu(E_*) < \mu(E) + \ep(d,r)/2 < \mu(E) +
\ep_-(d,r),$$
hence $\mu(E') \leq \mu(E)$ and $E$ is semistable.
\end{pf}

We thus get a morphism $\M_\al \longrightarrow \M_{r,d}$ which is the map of
Theorem \ref{thm:bh2} in case $(r,d)=1.$
By choosing the weights and quasi-parabolic structure correctly, we can fit
$\M_{r,d}$
and $\M_{r,d-1}$ into a chain diagram of maps as follows. Let $D=p$ and
$m=(1,\ldots, 1),$ and choose weights
$\al=(\al_1,\ldots,\al_r)$ with $\al_1 + \cdots + \al_r < \ep(r,d)/2.$ Suppose
$\al_1 < \eta < \al_2$ and
set $E'_*$ to be the parabolic bundle $E_*$ shifted by $\eta.$
Notice that $E'_*$ has degree $d-1,$ multiplicities
$m'=(1,\ldots,1),$ and weights $\al'=(\al_2 - \eta, \ldots, \al_r-\eta,
1-\eta+\al_1).$
If $\beta' \in V_{m'}$ is generic with
$\be'_1 + \cdots + \be'_r < \ep(r,d)/2,$ then we can connect
$\al'$ to $\be'$ in $V_{m'}$ by a line passing through a finite number of
hyperplanes $H_{\xi^1}, \ldots, H_{\xi^n},$ all of the form to which Theorem
\ref{thm:bh1} applies.
Choose weights $\al^i$ in the intermediate chambers and
$\ga^i \in H_{\xi^i}$ for $i=1,\ldots, n$ with $\al^n = \be'.$
Applying Theorem \ref{thm:bh1} each time we cross a hyperplane,
we get the following diagram:
\begin{equation} \label{eqn:hecke}
\begin{array}{lllllllll}\M_\al \cong \M_{\al'} & &&&\M_{\al^1} &  & &&
\!\!\!\M_{\be'}\\
 \psi \downarrow  &   \searrow  && \swarrow && \searrow && \swarrow &
\downarrow \psi' \!\!\!\! \\
\M_{r,d}  && \M_{\ga^1} &&& & \cdots & & \!\!\!\M_{r,d-1}
\end{array}
\end{equation}

\medskip \noindent
where, by the above proposition, the vertical maps $\psi$ and $\psi'$
have fibers the (full) flag variety over $\M^s_{r,d}$ and $\M^s_{r,d-1},$
respectively.
By Theorem \ref{thm:bh2}, $\psi$ is a
fibration which is locally trivial in
the Zariski topology provided $(r,d)=1,$ and the same follows for $\psi'$ if
$(r,d-1)=1.$

The second application of shifting is to extend the results of \cite{bh} to a
case which is
natural from the point of view of representations of Fuchsian groups but less
natural
from the point of view of parabolic bundles.
Assume for simplicity that $\mu(E_*) = 0$ and $D=p.$
Thus, $\deg E = -k$ for some $0 \leq k < r,$ and the relevant weight space
is
$$W_k = \{ (\al_1,\ldots,\al_r) \in \Delta^r \mid \al_1 + \cdots + \al_r = k
\}.$$
Consider the union
${\displaystyle \widetilde{W} = \bigcup_{k=0}^{r-1} W_k},$
where we identify
$$\partial_0 W_k = \{ \ga \in W_k \mid \ga_1=0 \}$$
with its companion set
$$\partial_1 \overline{W}_{k+1} = \{ \overline{\ga} \in \overline{W}_{k+1} \mid
\overline{\ga}_r=1 \}$$
via the identification
\begin{equation} \label{eqn:glue}
\partial_0 W_k \ni \ga = (0, \ga_2, \ldots, \ga_n) \sim (\ga_2, \ldots, \ga_n,
1) = \overline{\ga} \in \partial_1 \overline{W}_{k+1}.
\end{equation}
One can think of this set
$\widetilde{W}$ as the space of all weights
modulo shifting\footnote{Because every bundle can be shifted so that
$\mu(E_*)=0.$},
which in this case is just
the quotient $\SU(r)/Ad$ and which can be naturally
identified with the standard $r-1$ simplex.
{}From this point of view $\partial_0 W_k$ is
an interior hyperplane of $\widetilde{W}$
because it satisfies condition (\ref{eqn:hyper}).

However, Theorem \ref{thm:bh1} does not obviously carry over to this case
because points in $W_k$ and $W_{k+1}$ are weights on
parabolic bundles of different degrees.
Given a parabolic bundle of degree $-k,$
what is needed is a canonical procedure to construct a
parabolic bundle of degree $-(k+1).$
This is precisely what is provided by the shifting operation.
Thought of in terms of $ \widetilde{W},$
the following theorem extends Theorem \ref{thm:bh1} to
the case where $H_\xi= \partial_0 W.$

We use the notation $\M_\al(k,m)$ for the moduli space when
$E_*$ has degree $-k,$ multiplicities $m,$ and weights $\al.$

\begin{thm} \label{thm:bh3}
Suppose that $\ga \in \partial_0 W_k \cap V_m$ does not lie on any other
hyperplanes and that $\al \in W_k \cap V_m$ is a generic
weight near to $\ga.$ Choose $\eta \in {\Bbb R}$ with
$0<\eta < \ga_{m_1+1}.$ Define $\overline{\ga}\in \partial_1
\overline{W}_{k+1}$
as in \rom(\ref{eqn:glue}\rom).
Let $E'_*$ be $E_*$ shifted by $\eta,$ and denote the multiplicities of $E'_*$
by $m'.$
Set $k'= - \deg E' = k+m_1.$
Let  $\be \in W_{k'} \cap V_{m'}$ be generic near $\overline{\ga}.$
Then there are projective algebraic maps $\phi_\al$ and $\phi_\be$
$$\begin{array}{rcl}\M_\al(m,k)& & \M_\be(m',k')\\
\;\; \phi_\al \!\! \searrow \!\!\!\!\!\! && \!\!\!\!\! \swarrow \!\!
\phi_\be\\& \M_\ga(m,k)
\end{array}$$
satisfying the conclusion of Theorem \ref{thm:bh1}.
\end{thm}

\begin{pf}
By the choice of $\al, \be$ and $\eta,$ we see that $\al_{m_1} < \eta <
\al_{m_1+1},
\; \eta < \be_1$ and $\eta < \gamma_{m_1+1}.$
Consequently, the shifting operation defines the following isomorphisms:
\begin{eqnarray*}
\M_\al(m,k) &\cong& \M_{\al'}(m',k'),\\
\M_\be(m',k') &\cong& \M_{\be'}(m',k'),\\
\M_\ga(m,k) &\cong& \M_{\ga'}(m',k'),
\end{eqnarray*}
where $\al', \be', \ga' \in V_{m'}$ are defined from $\al, \be, \ga$ as in
(\ref{eqn:sft_wts}).
Now Theorem \ref{thm:bh1} applies to the shifted moduli spaces to prove the
theorem.
One can calculate $e_\al$ and $e_\be$ by applying formulas
(\ref{formula:ealpha}) and (\ref{formula:ebeta}) to $ \al', \be'$ and $\ga'.$
\end{pf}
\medskip
\noindent
{\it Remark.} Theorem \ref{thm:bh3} solves a problem mentioned at the end of
\cite{bh}
and extends the wall-crossing formula for knot invariants introduced in
\cite{boden3}.

\section{Rationality of moduli spaces of parabolic bundles}

Let $L$ be a holomorphic line bundle over a curve $X$ of genus
$g \geq 2.$ Denote by
\begin{enumerate}
\item[(i)] $\M_{r,L}$ the moduli space of semistable  bundles $E$ of rank $r$
with $\det E = L,$ and by
\item[(ii)] $\M_{\al, L}$ the moduli space of parabolic bundles $E_*$ with
weights $\al$ and
$\det E = L.$
\end{enumerate}
The main results of \S 4 hold for the moduli spaces with fixed determinant with
no essential difference. In view of
Theorem \ref{thm:bh2},
the goal is therefore to prove rationality with the coarsest possible choice of
flag structure.
At one extreme, we have the trivial flag, whose moduli space is exactly
$\M_{r,L}.$
Proposition 2 of \cite{newstead1} implies that $\M_{r,L}$ is rational if $\deg
L = \pm 1 \mod (r),$
and then
Theorem \ref{thm:bh2} and Proposition \ref{prop:qp} imply that $\M_{\al,L}$ is
also rational for any $\al \in V_m$
provided $\deg L = \pm 1 \mod (r).$

\begin{thm} \label{thm:rat1}
If $m(p)=(1,\ldots,1)$ for some $p \in D,$ then
$\M_{\al, L}$ is rational for all $\al \in V_m.$
\end{thm}
\begin{pf}
First, use Theorem \ref{thm:bh2} to reduce to the case $D=p$ by forgetting all
the other flag structures. If $E'_*$ denotes the bundle obtained by shifting
$E_*$
by some $\eta$ with $\al_1 < \eta < \al_2,$
then $\det E' = L' = L(-p).$ It follows that shifting by $\eta$ defines an
isomorphism from
$\M_{\al,L}$ to $\M_{\al',L'}.$
Repeated application of shifting puts us in the case $\deg L = 1 \mod(r),$ and
then Newstead's theorem and Theorem \ref{thm:bh2} imply that $\M_{\al,L}$ is
rational.
\end{pf}

The above argument works in slightly more generality.
We can always shift our bundle to be any of the $\E_x$ appearing in
the filtration (\ref{eqn:filtration1}) and illustrated in Figure 1.
Thus, whenever one of these terms in the filtration is of a degree to which
Newstead's theorem applies, the corresponding moduli space
of parabolic bundles is rational.

The next theorem is a considerable strengthening of the previous one.
\begin{thm} \label{thm:rat2}
If $m_i(p)=1$ for some $p \in D$ and some $1 \leq i \leq \kappa_p,$
then $\M_{\al,L}$ is rational for all $\al \in V_m.$
\end{thm}
Before delving into the proof of this theorem, we mention some interesting
consequences.
Recall first the following definition.
\begin{defn} \label{defn:stablerat}
A variety $V$ is  stably rational of level $k$ if $V \times {\Bbb P}^k$ is
rational.
The level is the smallest integer $k$ with this property.
\end{defn}

The following result, with a weaker bound on the level, was proved in
\cite{ballico}.
\begin{cor} \label{cor}
For $(r,d)=1, \; \M_{r,L}$ is stably rational with level $k \leq r-1.$
\end{cor}
\begin{pf} Theorem  \ref{thm:rat2} implies that $\M_{\al,L}$ is rational, where
$m(p)=(r-1,1),$ and
Theorem \ref{thm:bh2} shows that
$\M_{\al,L}$ is birational to $\M_{r,L} \times {\Bbb P}^{r-1},$  which proves
the corollary.
\end{pf}

We now apply this last result to Conjecture \ref{conj}.
\begin{cor} \label{cor:rat}
Suppose $(r,d)=1.$ By tensoring with a line bundle, we can assume that $0 < d <
r.$
If either $(g,d) =1$ or $(g,r-d)=1,$ then $\M_{r,L}$ is rational.
\end{cor}
\begin{pf}
Suppose first that $(g,r-d)=1.$ Let $L$ be a line bundle of degree $r(g-1)+d.$
Then Newstead's construction applies and proves that
$\M_{r,L}$ is birational to $\M_{r-d,L} \times {\Bbb P}^\chi,$ where $\chi =
(g-1)(r^2-(r-d)^2).$
But the above corollary implies that $\M_{r-d,L}$ is stably rational with level
$k \leq r-d-1 \leq \chi,$
hence $\M_{r,L}$ is rational.

The case $(g,d)=1$ follows by the same argument after applying duality, which
interchanges
$(r,d)$ and $(r,r-d).$
\end{pf}
{\it Remark.} Conjecture \ref{conj} was previously known \cite{newstead1} in
the following three cases:
\begin{enumerate}
\item[(i)] $d = \pm 1 \mod (r),$
\item[(ii)] $(r,d)=1$ and $g$ a prime power, and
\item[(iii)] $(r,d)=1$ and the two smallest distinct primes factors of $g$ have
sum
greater than $r.$
\end{enumerate}
Corollary \ref{cor:rat}
applies in each case. More importantly, it applies
in many cases not covered by (i), (ii) or (iii).
In fact, for a given $r$ and $d$ with $(r,d)=1,$ one can easily
list those $g$ for which the conjecture remains open.
For example, if $r= 110$ and $d=43,$ then Corollary \ref{cor:rat} applies
as long as $g$ is not a multiple of $d \cdot (r-d)= 43 \cdot 67=2881.$

 \medskip
\noindent
{\it Proof of Theorem.}
Set $d= \deg L.$
The theorem is clearly true for $r=1$ and follows
from Theorem \ref{thm:rat1} for $r=2,$ so assume $r>2.$
Notice that by tensoring with a line bundle, we can suppose $$r(g-1) < d \leq
rg.$$
By Theorem \ref{thm:bh2}, we can again assume that
$D=p,$ and by shifting and another application of Theorem \ref{thm:bh2}, if
necessary,
we can arrange it so that $m(p)=(r-1,1).$
Write $$\al = \al(p) = (\overbrace{\al_1, \cdots, \al_1}^{r-1}, \al_2).$$
Proposition \ref{prop:fine} implies that $V_m$ contains a generic weight
and that $\M_{\al,L}$ parametrizes  a universal family $\U_*^\al.$
By Proposition \ref{prop:qp}, the birational type of $\M_{\al,L}$ is
independent of
choice of compatible weights, so
we can assume that the weights are small enough to satisfy the hypothesis of
Proposition \ref{prop:num}
(this comes up at various technical points in the argument, e.g. the proof of
Claim \ref{claim}).

\medskip\noindent
Consider the following two cases.

\medskip\noindent
{\sc Case I:} \quad $d=rg.$ \quad
Choose $\eta$ with $\al_1 < \eta < \al_2,$ and let $E'_* = E_*[\eta]_*.$
Denote the weights of $E'_*$ by $\al'$ as in (\ref{eqn:sft_wts}).
If $\det E = L,$ then $\det E' = L' = L(-(r-1)p)$ has degree $d'= d-(r-1).$
Since $d' = 1 \mod(r),$  Proposition 2 of \cite{newstead1}
implies that $\M_{r,L'}$ is rational, and Theorem \ref{thm:bh2} then implies
that
$\M_{\al'\!,L'}$ is also rational.
Rationality of $\M_{\al,L}$ now follows from the isomorphism of the
moduli spaces $\M_{\al,L} \cong \M_{\al'\!,L'}$ defined by shifting by $\eta.$


\medskip\noindent
{\sc Case II:} \quad $r(g-1) < d < rg.$ \quad
The idea is to use induction to construct a nonempty, Zariski-open subset $\M$
of affine space of dimension
$(r^2-1)(g-1)+r-1 \;  (= \dim \M_{\al,L})$
and a family of stable parabolic bundles $\U_*$ parametrized by $\M$ with $\det
\U_{\xi,*} = L$ for all
$\xi \in \M.$
The universal property of $\U_*^\al$ then gives a map $\psi_{\U_*} : \M
\longrightarrow \M_{\al,L}.$
If, in addition, we have $\U_{\xi_1,*} \cong \U_{\xi_2,*}\; \Leftrightarrow \;
\xi_1 = \xi_2,$ then $\psi_{\U_*}$ is injective
and rationality of $\M_{\al,L}$ follows from that of $\M$ and the dimension
condition.

Set $r' = rg-d, \, r'' = r-r'$ and $\al'= (\overbrace{\al_1, \cdots,
\al_1}^{r'-1}, \al_2).$
Assume that both $\al$ and $\al'$ are generic.
Let $\U_*^{\al'}$ be the universal family parametrized by $\M_{\al'\!,L}$
and $I_*= {\cal O}_X[\al_1]_*$ be the trivial parabolic line bundle with weight
$\al_1.$
If $e' = [E'_*] \in \M_{\al'\!,L},$
then because ${E'_*}^\vee\otimes I_*$ is a stable parabolic bundle of negative
parabolic degree,
$h^0({E'_*}^\vee \otimes I_*)=0$ and
\begin{equation} \label{eq:rank}
 n \stackrel{\hbox{\scriptsize\rm def}}{=} h^1({E'_*}^\vee \otimes I_*) =
(2r'+r'')(g-1) + r'' + 1
\end{equation}
is independent of $e'.$
Since $\U^{\al'}_{e',*} \cong E'_*,$ it follows that
$$(R^1 \pi_{\M_{\al'\!,L}}) (({\U_*^{\al'}})^\vee \otimes \pi_X^*(I_*))$$
is locally free. The associated vector bundle $V
\stackrel{\pi}{\longrightarrow} \M_{\al'\!,L}$
has rank $n$ and fiber over $e'$ naturally isomorphic to $H^1({E'_*}^\vee
\otimes I_*).$

Let $\U_*'=(\pi^{r''} \times 1_X)^* (\U_*^{\al'})$ be the pullback family
and $\I^{\oplus r''}_* =\pi_X^* I_*^{\oplus r''}$ the trivial family,
where $\pi^{r''}:V^{\oplus r''} \longrightarrow \M_{al'\!,L}$.
There is an extension
\begin{equation} \label{eqn:universal_extension}
0 \longrightarrow \I_*^{\oplus r''} \longrightarrow \U_* \longrightarrow \U_*'
\longrightarrow 0
\end{equation}
of families over $V^{\oplus r''} \times X,$
such that, for $\xi \in V^{\oplus r''}_{e'},$
$\U_{\xi,*}$ is the parabolic bundle $E^\xi_*$
described as the short exact sequence
\begin{equation} \label{eq:ses}
0 \longrightarrow I^{\oplus r''}_*  \longrightarrow E^\xi_* \longrightarrow
E'_* \longrightarrow 0
\end{equation}
corresponding to the extension class $\xi \in H^1({E'_*}^\vee \otimes I^{\oplus
r''}_*).$

Using stability of $E'_*$ and triviality of $I^{\oplus r''}_*,$
it follows that
$$\Aut(E'_*) \times \Aut(I^{\oplus r''}_*) \cong {\Bbb C}^* \times \GL(r'',
{\Bbb C}).$$
This group acts naturally as fiber-preserving maps on the bundle $V^{\oplus
r''}$  since
$$V^{\oplus r''}_{e'} \cong H^1({E'_*}^\vee \otimes I^{\oplus r''}_*) =
H^1({E'_*}^\vee \otimes I_*)^{\oplus r''},$$
and two extension classes $\xi_1$ and $\xi_2$ in the same orbit
have associated bundles $E^{\xi_1}$ and $E^{\xi_2}$ which are isomorphic.
We can ignore the ${\Bbb C}^*$ action here
because $(z,1)\cdot \xi = (1,z) \cdot \xi$ for $z \in {\Bbb C}^*$ and $\xi \in
V^{\oplus r''}.$

Using the inductive
hypothesis and local triviality of $V,$ we can choose
a nonempty Zariski-open subset $\M'$ of $\M_{\al'\!,L}$ isomorphic to a
Zariski-open subset of affine space of dimension $({r'}^2-1)(g-1)+r'-1$ 
such that $V|_{\M'} \cong \M' \times H^1({E'_*}^\vee \otimes I_*)$ ($E'_*$ is
fixed).
Lemma 2 of \cite{newstead1} applies here and produces a Zariski-open subspace
$\M'\times W$ of
$V^{\oplus r''}|_{\M'}$ invariant under the group action,
and an affine subspace $A \subset W$
so that every orbit in $W$ intersects $A$ precisely once.
In fact, $A$ can be chosen as a Zariski-open subset of the Grassmannian
$G(r'',n).$
In any case, it should be clear that $A$ has dimension $r''(n-r'').$
Using equation (\ref{eq:rank}) and the fact that $r'+r'' = r,$ we see
that $\M' \times A$ is a Zariski-open subset of affine space of dimension
\begin{eqnarray*}
\dim \M' \times A &=& ({r'}^2-1)(g-1)+r'-1 + r''(n - r'')\\
&=& ({r'}^2-1)(g-1)+r'-1 + r''((2r'+r'')(g-1) + 1) \\
&=& (r^2-1)(g-1)+r-1.
\end{eqnarray*}

Let $\M$ be the subset of $V^{\oplus r''}$ defined by
$$\M = \{ \xi \in \M' \times A \mid H^1(\U_{\xi,*}) = 0 \},$$
and consider the bundle $\U_*$ restricted to $\M,$ which we continue to denote
$\U_*.$
For $\xi \in V^{\oplus r''},$ let $E^\xi_* = \U_{\xi,*}.$
Clearly $E^\xi_*$ is a parabolic bundle with weights $\al$ and determinant $L,$
thus $\M$ parametrizes a family of parabolic bundles.
By the upper semi-continuity theorem, $\M$ is Zariski-open in $\M' \times A.$

We claim that $\M$ is nonempty.
Fix $e'=[E'_*] \in \M'$ and consider the set
$$N= \{ \xi \in H^1({E'_*}^\vee \otimes I^{\oplus r''}_*) \mid h^1(E^\xi_*)=0
\}.$$
If $N \cap A \neq \emptyset,$ then $\M$ is nonempty. Clearly, $N$ is invariant
under the
action of $\GL(r'',{\Bbb C}),$ so it is enough to show $N \cap W \neq
\emptyset.$
There is a natural map
$$\delta : H^1({E'_*}^\vee \otimes I^{\oplus r''}_*) \times H^0(E'_*)
\longrightarrow H^1(I^{\oplus r''}_*)$$
with $\delta_\xi = \delta(\xi,\cdot) : H^0(E'_*) \longrightarrow H^1(I^{\oplus
r''}_*)$ the
coboundary map of the long exact sequence in homology of
(\ref{eq:ses}).
Now $H^0(E'_*)=H^0(E'),$ and since $\al_1+(r'-1) \al_2 < \ep(r,d)/2,$ by
Proposition \ref{prop:num},
$E'$ is semistable as a non-parabolic bundle. Serre duality implies that
$h^1(E')=h^0({E'}^\vee \otimes K),$ and we compute
\begin{eqnarray*}
 \deg({E'}^\vee \otimes K) &=& -d +r'(1-g) \\
&\leq& (r+r')(1-g) - r'',
\end{eqnarray*}
which is negative since $r'' \geq 1$ and $g \geq 2.$ This implies that
$h^1(E'_*) =0,$
and Riemann-Roch implies that $h^0(E'_*) = r''g.$
Because $h^1(I^{\oplus r''}_*)=r''g,$
we see that
$$\xi \in N \Longleftrightarrow H^1(E^\xi_*) =0 \Longleftrightarrow \delta_\xi
\hbox{ is an isomorphism.}$$
But $\delta$ is obviously onto and
$\dim (\ker \delta) = r''n.$
The set $N$ has complement
$$N^c = \{ \xi \in H^1({E'_*}^\vee \otimes I^{\oplus r''}_*)
\mid \delta(\xi,s)=0 \hbox{ for some }0 \neq s \in H^0(I^{\oplus r''}_*) \}.$$
But $\delta(\xi,s) = 0 \Rightarrow \delta(\xi, zs)=0$ for all $z \in {\Bbb C},$
which shows that
the map $\ker \delta \longrightarrow N^c$ has fibers of dimension $\geq 1.$
Hence $\dim N^c \leq \dim(\ker \delta) -1 < r''n,$ and we see that $N$ is
nonempty and Zariski-open.
Thus $N \cap W \neq \emptyset$ and it follows that $\M$ is nonempty.

We now prove that $\M$ parametrizes a family of stable parabolic bundles,
using again the inequality
$(r-1) \al_1 + \al_2 < \ep(r,d)/2$ and Proposition \ref{prop:num}.
\begin{claim} \label{claim}
\begin{enumerate}
\item[(i)] $E^\xi_*$ is stable for all $\xi \in \M.$
\item[(ii)] $E^{\xi_1}_* \cong E^{\xi_2}_* \Longleftrightarrow
\GL(r'', {\Bbb C}) \cdot \xi_1 = \GL(r'', {\Bbb C}) \cdot \xi_2$ for all
$\xi_1, \xi_2 \in \M.$
\end{enumerate}
\end{claim}
\begin{pf}
(i)  Suppose to the contrary that $E^\xi_*$ is not parabolic stable for some
$\xi \in \M.$
Let $G_*$ be a rank $s$ parabolic subbundle of $E^\xi_*$ with $\mu(G_*) >
\mu(E^\xi_*).$
Then $\mu(G) \geq \mu(E^\xi),$ since otherwise
$$\mu(G_*) < \mu(G) + \ep(d,r)/2 < \mu(E^\xi) < \mu(E^\xi_*).$$
As in the argument of Lemma 6 of Newstead, the map
$G \longrightarrow E'$ has a factorization as
$G \rightarrow G^1 \rightarrow G^2 \rightarrow E'$  and the arguments there
give the following inequalities:
\begin{eqnarray} \label{eq:1}
\deg(G^2) &\geq& \deg(G) \geq \frac{sd}{r}, \\
\label{eq:2}
\rank(G^2) &\leq& \rank(G) - h^0(G) \leq \frac{sr'}{r}.
\end{eqnarray}
These imply that $\mu(G^2) - \mu(E') \geq 0.$
But $E'_*$ is parabolic stable, so by Proposition \ref{prop:num},
$E'$ is semistable and $\mu(G^2)=\mu(E').$
Thus, we must have equalities
in equations (\ref{eq:1}) and (\ref{eq:2}), in particular $\mu(G) =
\mu(E^\xi)$.
But since $\mu(G_*) > \mu(E^\xi_*),$ we see that $G_*$ must inherit the weight
$\al_2,$
which implies that $G^2_*$ also inherits $\al_2,$ and it now follows that
$$\mu(G^2_*) - \mu(E'_*) = \frac{(s_2 -1) \al_1 + \al_2}{s_2} - \frac{(r'-1)
\al_1 + \al_2}{r'} > 0,$$
where $s_2 = \rank G^2 < r'.$ This contradicts the parabolic stability of
$E'_*$ and completes the proof of part (i).

\medskip\noindent
(ii) Since $\Leftarrow$ is true independent of the vanishing of $H^1,$ we only
prove $\Rightarrow.$
Suppose $E^{\xi_1}_* \cong E^{\xi_2}_*$ and set $\pi_X(E^{\xi_i}_*) =
e_i'=[{E^i_*}'] \in \M_{\al'\!,L}.$
Notice that $h^1(E^{\xi_i}_*) = 0,$
and so $h^0(E^{\xi_i}_*) = \chi(E^{\xi_i}_*) = r''.$ It follows that
every holomorphic section of $E^{\xi_i}_*$ has its image contained in
$I^{\oplus r''}_*.$
Hence any isomorphism $\varphi : E^{\xi_1}_* \longrightarrow E^{\xi_2}_*$
defines
a commutative diagram
$$\CD
0 @>>> I^{\oplus r''}_*  @>>>   E^{\xi_1}_*  @>>>  {E^1_*}' @>>> 0 \\
   &&   @VV{\varphi''}V         @VV{\varphi}V  @VV{\varphi'}V \\
0 @>>> I^{\oplus r''}_*  @>>>   E^{\xi_2}_*  @>>>  {E^2_*}' @>>> 0 \\
\endCD $$
where both $\varphi'$ and $\varphi''$ are isomorphisms, and
so $\xi_2 = (\varphi' \times \varphi'') \cdot \xi_1.$
\end{pf}

Part (i) of the claim and the universal property of $\U_*^\al$
gives a map $\M \stackrel{\psi}{\longrightarrow} \M_{\al,L},$ which is
injective by part (ii).
Since $\M$ is nonempty,
$\dim \M = \dim  \M_{\al,L},$ so rationality of $\M_{\al,L}$ follows from that
of $\M.$
This concludes the proof in Case II.
$\quad \Box$

\medskip
\noindent
{\it Remark.}
We had originally hoped to prove rationality of $\M_{\al,L}$ with the weaker
hypothesis that
$\al$ is generic,
but the argument does not hold in this generality.
For consider the case $D=p.$ By tensoring with a line bundle and shifting, we
can assume that
$$r(g-1) < d \leq r(g-1) + m_1.$$ Hence, the subbundle split off in the
induction
is again a sum of parabolic line bundles with the same weights.
The difficulty is in proving that
the quotient $E'_*$ has {\it generic} weights $\al'.$

Proposition \ref{prop:fine} implies that
$E'_*$ admits a generic weight if and only if the elements of the set $\{d,
m'_i(p) \}$ greatest common divisor equal to one.
The statement
$$(d,m_1,\ldots,m_\kappa)=1 \Rightarrow (d,m'_1,\ldots,m'_\kappa)=1,$$
which is what we would need to prove here, is unfortunately false
(notice that $m'_1=m_1-d+r(g-1)$ and $m'_i=m_i$ otherwise).

\medskip\noindent
{\it Acknowledgements.} Both authors would like to express their gratitude to
the
Max-Planck-Institut f\"ur Mathematik for providing financial support. The first
author
is also grateful to the Institut des Hautes \'Etudies
Scientifique for partial support.
We would also like to thank I.\ Dolgachev and L.\ G\"ottsche
for helpful discussions.

\end{document}